\documentclass[12pt]{article}
\usepackage{amssymb,amsmath,latexsym}
\usepackage{xcolor}
\def\ZZ{{\mathbb Z}}
\def\QQ{{\mathbb Q}}
\def\RR{{\mathbb R}}
\def\CC{{\mathbb C}}

\newtheorem{formula}{}[section]
\newtheorem{proposition}[formula]{Proposition}
\newtheorem{definition}[formula]{\indent Definition}
\newtheorem{corollary}[formula]{\indent Corollary}
\newtheorem{remark}[formula]{\indent Remark}
\newtheorem{lemma}[formula]{\indent Lemma}
\newtheorem{theorem}[formula]{\indent Theorem}

\def\thrm{\begin{theorem}}
\def\thrml#1{\begin{theorem}\label{#1}}
\def\ethrm{\end{theorem}}
\def\rmrk{\begin{remark}}
\def\rmrkl#1{\begin{remark}\label{#1}}
\def\ermrk{\end{remark}}
\def\dfntn{\begin{definition}}
\def\dfntnl#1{\begin{definition}\label{#1}}
\def\edfntn{\end{definition}}
\def\nmrt{\begin{enumerate}}
\def\enmrt{\end{enumerate}}

\def\qtn{\begin{equation}}
\def\qtnl#1{\begin{equation}\label{#1}}
\def\eqtn{\end{equation}}
\def\lmm{\begin{lemma}}
\def\lmml#1{\begin{lemma}\label{#1}}
\def\elmm{\end{lemma}}
\def\crllr{\begin{corollary}}
\def\crllrl#1{\begin{corollary}\label{#1}}
\def\ecrllr{\end{corollary}}

\begin{document}
\title{}
\date{}
\maketitle
\vspace{-0,1cm} \centerline{\bf Degree of 1-dimensional tropical prevariety}
\vspace{7mm}
\author{
\centerline{Dima Grigoriev}
\vspace{3mm}
\centerline{CNRS, Math\'ematique, Universit\'e de Lille, Villeneuve
d'Ascq, 59655, France} \vspace{1mm} \centerline{e-mail:\
dmitry.grigoryev@univ-lille.fr } \vspace{1mm}
\centerline{URL:\ http://en.wikipedia.org/wiki/Dima\_Grigoriev} }

\begin{abstract}
For a tropical prevariety $V\subset \RR^n$ (being a finite union of rational polyhedra) we define a tropical Hilbert function $TH_V(k)$ to be the maximal number of tropically independent on $V$ among tropical monomials with degrees at most $k$. In case $\dim V=1$ we define the tropical degree as
$$degT (V):=\lim_{k\to \infty} \frac{TH_V(k)}{k}$$
\noindent and prove existence of the limit. We calculate explicitly (a modification of) the tropical degree of $V$ when $V=\cup_l V_l$ is a star, i.e. the union of rays $V_l$ with a common apex.

\end{abstract}

{\bf keywords}: min-plus prevariety, tropical Hilbert function, tropical degree


\section*{Introduction}

Recall that a tropical (or min-plus) polynomial $f=\min_{1\le i\le s}\{a_{i,1}x_1+\cdots + a_{i,n}x_n+a_{i,0}\}$ where $0\le a_{i,j}\in \ZZ \cup \{\infty\},\ 1\le i\le s, 1\le j\le n,\ a_{i,0}\in \RR \cup \{\infty\}$ \cite{MS}, \cite{J}. (The linear functions) $a_{i,1}x_1+\cdots + a_{i,n}x_n+a_{i,0},\ 1\le i\le s$ are called tropical monomials, $a_{i,1}+\cdots+ a_{i,n}$ are called their tropical degrees, and the maximum of the latter over $1\le i\le s$ is called the tropical degree of $f$.

A {\it min-plus prevariety} $V:=V(f_1=g_1,\dots,f_m=g_m)\subset \RR^n$ of a system of min-plus equations $f_1=g_1,\dots,f_m=g_m$ (where $f_1,\dots,f_m,g_1,\dots,g_m$ are tropical polynomials) is the set of all  points $x\in \RR^n$ satisfying $f_1(x)=g_1(x),\dots, f_m(x)=g_m(x)$. The category of min-plus prevarieties coincides with the category of tropical prevarieties (an explicit transition between these two categories one can find  e.g. in \cite{GP18}), and also coincides with the category of unions of rational polyhedra (so, polyhedra determined by linear inequalities with rational coefficients) \cite{MS}.

For a min-plus prevariety $V$ the paper \cite{BE} introduces {\it congruence classes} of tropical polynomials such that tropical polynomials belong to the same class iff they coincide on $V$. The family of congruence classes is denoted by $E(V)$. An explicit description of $E(V)$ in terms of $f_1,\dots,f_m,g_1,\dots,g_m$ was conjectured in \cite{BE} and proved in \cite{JM}. One can treat this description as a min-plus version of Hilbert's strong Nullstellensatz. We mention that in \cite{GP18} a tropical version of Hilbert's weak Nullstellensatz is provided (further developments one can find in \cite{MR}, \cite{ABG}). Still, it is not clear whether there is a tropical version of Hilbert's strong Nullstellensatz.

Briefly remind that for an algebraic variety $U\subset \CC^n$ its Hilbert function $H_U(k)$ is defined as the dimension of the vector space of polynomials on $U$ with degrees at most $k$. In other words, it equals $\dim (\CC[x_1,\dots,x_n]_k/I(U)_k)$, where $I(U)$ is the ideal of all polynomials vanishing on $U$. It is well known that $H_U(k)$ is a polynomial for sufficiently big $k$, with degree equal to $\dim U$ and the leading coefficient equal to the degree of $U$.

In \cite{MR} tropical ideals in the semiring of tropical polynomials are introduced. They have some common features with ideals in the polynomial ring: in particular, the Hilbert function is defined for tropical ideals, and it is proved that the Hilbert function for a homogeneous tropical ideal is eventually a polynomial. Moreover, for a tropicalization of an algebraic variety $U$ (say, over the field of Newton-Puiseux series) its tropical Hilbert function coincides with Hilbert function of $U$.

In the present paper we define a {\it tropical Hilbert's function} $TH_V(k)$ as the maximal number of tropically independent (see Definition~\ref{dependent}) on a tropical prevariety $V$ among the tropical monomials of degree at most $k$ (see Definition~\ref{tropical_hilbert}). Unfortunately, the notion of tropical independence does not fulfill the matroid property (cf. \cite{DSS}). However, in case $\dim V=1$ we define the {\it tropical degree}
$$degT (V):= \lim_{k\to \infty} \frac{TH_V(k)}{k}$$
\noindent and prove existence of the limit (Theorem~\ref{degree}). We conjecture that $TH_V(k)$ is a polynomial (of degree $\dim V$) for sufficiently large $k$, in this case (when $\dim V=1$) $degT (V)$  is the leading coefficient of this (linear) polynomial.

We define a modification $degT^{(\Box)} (V)$ of the tropical degree where in the definition of $TH_V^{(\Box)} (k)$ the maximum is taken over the number of tropically independent on $V$ monomials among $\{i_1x_1+\cdots +i_nx_n\ |\ 0\le i_1,\dots, i_n\le k\}$. 
We calculate explicitly $degT^{(\Box)} (V)$ when $V=\cup_l V_l\subset \RR^2$ is a "star", i.e. each $V_l$ is a 
ray spanned by a vector $(p_l,q_l)\in \ZZ^2$ with relatively prime $p_l,q_l$, in addition, all $V_l$ have a common apex. We show (Theorem~\ref{star}) that
$$degT^{(\Box)} (V)= \max\{\sum_{p_l>0} p_l,\ -\sum_{p_l<0} p_l\} + \max\{\sum_{q_l>0} q_l,\ -\sum_{q_l<0} q_l\}.$$

{\bf Questions.} 1) Is there a geometric meaning of the tropical degree?

2) Is the degree always an integer?

3) For $\dim V=m$ does it exist the limit
$$\lim_{k\to \infty} \frac{TH_V(k)}{k^m}?$$


We mention that in \cite{EG} a different version of a tropical Hilbert's function was suggested for univariate tropical polynomials. It was proved in \cite{EG} that this tropical Hilbert's function coincides with the floor function of a linear polynomial with rational coefficients (for sufficiently big argument). The leading coefficient of this linear function is called the tropical entropy, in \cite{G20} a criterion is provided when the entropy equals zero. In \cite{G24} the entropy is studied for multivariate tropical polynomials. 

\section{Tropical Hilbert function of a tropical prevariety}

\begin{definition}\label{dependent}
Let $V\subset \RR^n$ be a min-plus prevariety. We say that tropical polynomials $f_1,\dots, f_m$ in $n$ variables $X_1,\dots, X_n$ are {\it tropically independent on $V$} if there there exist $b_1,\dots, b_m\in \RR$ and $v_1,\dots,v_m\in V$ such that
\begin{equation}\label{3}
b_{j_0}+f_{j_0} (v_{j_0})< b_j+f_j(v_{j_0}),\ 1\le j\neq j_0\le m.
\end{equation}
In this case we call the tropical polynomial $\min_{1\le j\le m} \{b_j+f_j\}$ a {\it certificate of independence of $f_1,\dots,f_m$}.
\end{definition}

One can verify that this concept of independence is equivalent to the following one suggested in \cite{GG}. Namely, $f_1,\dots, f_m$ are independent if a certain tropical linear combination of $f_1,\dots, f_m$ cannot be represented as a tropical linear combination of a proper subset of $f_1,\dots, f_m$. The equivalence follows from the observation that for a certificate it holds that if $\min_{1\le j\le m} \{b_j+f_j\}=\min_{1\le j\le m} \{b_j'+f_j\}, b_j'\in \RR \cup \{\infty\}$ then $b_j=b_j', 1\le j\le m$.

\begin{definition}\label{tropical_hilbert} For a tropical prevariety $V\subset \RR^n$ denote by the tropical Hilbert function $TH_V(k)$  the maximal number of independent on $V$ tropical monomials among the set
$$M_k:=\{i_1x_1+\cdots + i_nx_n\ |\ 0\le i_1+\dots + i_n \le k\}.$$
\end{definition}

\begin{remark}
For a min-plus prevariety $V\subset \RR^n$ define a congruence $E(V)$ on $n$-variate tropical polynomials as follows: two tropical polynomials $f,g$ are congruent iff $f|_V=g|_V$ \cite{BE}, \cite{JM}. One can identify the space $T_n(k)$ of tropical polynomials being tropical linear combinations of tropical monomials from $M_k$ with $\RR^{n+k \choose n}$. 
One can verify
that the dimension of the family of congruence classes  $T_n(k) \cap E(V)$ coincides with $TH_V(k)$ (cf. Theorem 4.2 \cite{DSS}). Indeed, this dimension does not exceed $TH_V(k)$ since for any tropical polynomial $f\in T_n(k)$ there exists a tropical polynomial
$$g:=\min_{1\le j\le r} \{b_j+m_j\ |\ m_j\in M_k, b_j\in \RR\}\in T_n(k)$$
\noindent such that $f|_V=g|_V$ for suitable $r\le TH_V(k)$. 
\end{remark}  

One can verify the following lemma with the help of the duality of linear programming (cf. Proposition 8.3 \cite{DSS}).

\begin{lemma}\label{diagonal}
Let $s\times s$ matrix $(w_{i,j})_{1\le i,j\le s}$ with $w_{i,j}\in \RR$ be tropically non-singular with the sum of the diagonal entries $w_{1,1}+\cdots + w_{s,s}$ being the strict minimum among all $s!$ matchings. Then there exist $w_1,\dots,w_s\in \RR$ such that
$$w_{i,i}+w_i<w_{l,i}+w_l,\ 1\le i\neq l\le s.$$ 
\end{lemma}

Now we rephrase the concept of independency in terms of the tropical rank of a suitable matrix.  For points $v_1,\dots, v_s\in \RR^n$ and tropical polynomials $f_1,\dots, f_m$ in $n$ variables consider $m\times s$ matrix $A(f_1,\dots, f_m; v_1,\dots, v_s)$ whose $(i, j)$-entry equals $f_i(v_j),\ 1\le i\le m, 1\le j\le s$. 

\begin{proposition}\label{trank}
The tropical polynomials $f_1,\dots, f_m$ are independent on a tropical variety $V\subset \RR^n$ iff there exist points  $v_1,\dots, v_m\in V$ such that the tropical rank 
\begin{eqnarray}\label{11}
trk A(f_1,\dots, f_m; v_1,\dots, v_m)=m.
\end{eqnarray}
Therefore, the maximal number of independent among $f_1,\dots, f_m$ on $V$ equals the maximal tropical rank $trk A(f_1,\dots, f_m; v_1,\dots, v_s)$ over all families of points $v_1,\dots, v_s \in V$.
\end{proposition}

{\bf Proof}. First assume that $f_1,\dots, f_m$ are independent on $V$, and $v_1,\dots, v_m \in V, b_1,\dots, b_m\in \RR$ are such that 
\begin{eqnarray}\label{10}
b_i+f_i(v_i)< b_j+f_j(v_i),\ 1\le i\neq j\le m.
\end{eqnarray}
Hence (\ref{11}) holds.
For a permutation $p\in Sym(m)$ denote 
\begin{eqnarray}\label{12}
A(f_1,\dots, f_m; v_1,\dots, v_m)_p:= f_1(v_{p(1)})+\cdots + f_m(v_{p(m)}).
\end{eqnarray}

Conversely, let (\ref{11}) be true. Asuume w.l.o.g. that 
$$A(f_1,\dots, f_m; v_1,\dots, v_m)_e:=f_1(v_1)+\cdots +f_m(v_m)< A(f_1,\dots, f_m; v_1,\dots, v_m)_p$$
\noindent for any non-unit permutation $p\neq e \in Sym(m)$. Lemma~\ref{diagonal} implies the existence of $b_1,\dots, b_m$ satisfying (\ref{10}). $\Box$ \vspace{2mm}





\section{A tropical Hilbert function of a polyhedron}


Throughout this section we assume that a min-plus prevariety $V\subset \RR^n$ is a rational polyhedron of dimension $m$. Denote by $L\subset \RR^n$ the parallel shift of the linear hull of $V$ such that $0\in L$. Note that two tropical monomials $\langle a,x \rangle,\langle b,x \rangle$ (where $a,b \in \ZZ_{\ge 0}^n$) coincide on $V$ up to an additive constant iff $(a-b) \perp L$. In the latter case we say that vectors $a,b$ (or tropical monomials $\langle a,x \rangle,\langle b,x \rangle$) are equivalent with respect to $L$. The number of equivalence classes of tropical monomials in $M_k$ we denote by $t(k):= t_L(k)$.

We need two auxiliary lemmas. Similar to \cite{GP20} we say that two families of real vectors $\{u_i=(u_{i,1},\dots,u_{i,m}),\ 1\le i\le s\}$ and  $\{v_i=(v_{i,1},\dots,v_{i,m}),\ 1\le i\le s\}$ are {\it co-ordered} if the inequality $u_{i_1,j}\le u_{i_2,j}$ is equivalent to $v_{i_1,j}\ge v_{i_2,j}$ for all $1\le i_1,i_2 \le s, 1\le j\le m$.

\begin{lemma}\label{co-ordered} (cf. Theorem 2 \cite{GP20}).
Let a pair of families of pairwise distinct vectors $\{u_i, 1\le i\le s\}$ and $\{v_i, 1\le i\le s\}$ be co-ordered. Then for any non-identic permutation $\pi \in Sym(s)$ it holds
\begin{equation}\label{2}
\sum_{1\le i\le s} \langle u_i, v_i \rangle < \sum_{1\le i\le s} \langle u_i, v_{\pi(i)} \rangle.
\end{equation}
\end{lemma}

{\bf Proof of the lemma}. Denote $S_{j,\pi}:= \sum_{1\le i\le s} u_{i,j}v_{\pi(i),j},\ 1\le j\le m$, then the right-hand side of (\ref{2}) equals $\sum_{1\le j\le m} S_{j,\pi}$. For each fixed $1\le j\le m$ the sum $S_{j,\pi}$ does not decrease with augmenting the set of inversions of $\pi$, i.e. the pairs $1\le i_1 < i_2\le s$ for which holds $\pi(i_1)> \pi(i_2)$ (the number of inversions is denoted by $l(\pi)$). Indeed, w.l.o.g. one can assume that $u_{1,j}\ge u_{2,j}\ge \cdots \ge u_{s,j}$ (thereby, $v_{1,j}\le v_{2,j}\le \cdots \le v_{s,j}$). One can represent $\pi$ as a composition of the minimal number (equal to $l(\pi)$) of transpositions of the form $(r,r+1),\ 1\le r<s$. If $l(\pi \circ (r,r+1))> l(\pi)$ then
$$S_{j,\pi \circ (r,r+1)}-S_{j,\pi}=u_{r,j}v_{\pi(r+1),j}+u_{r+1,j}v_{r,j}-u_{r,j}v_{\pi(r),j}-u_{r+1,j}v_{\pi(r+1,j}\ge 0.$$

This implies the non-strict inequality in (\ref{2}). Moreover, this argument justifies the following statement. Let $p_1\ge \cdots \ge p_s,\ q_1\le \cdots \le q_s$ be two sequences of the reals. Then 
$\sum_{1\le i\le s} p_iq_i \le \sum_{1\le i\le s} p_iq_{\pi(i)}$ for any  permutation $\pi$. When, in addition, $p_{i_1}=p_{i_2}$ is equivalent to that $q_{i_1}=q_{i_2}$ for any $1\le i_1, i_2\le s$, it holds $\sum_{1\le i\le s} p_iq_i < \sum_{1\le i\le s} p_iq_{\pi(i)}$, unless the partition of the set $\{1,\dots,s\}$ into the subsets in which $i_1,i_2$ belong to the same subset iff $p_{i_1}=p_{i_2}$, is invariant under $\pi$.

Since the vectors $\{u_i, 1\le i\le s\}$ are pairwise distinct, we conclude that for any non-identical permutation $\pi$ there exists $1\le j\le m$ such that $\sum_{1\le i\le s} u_{i,j}v_{i,j} < S_{j,\pi}$ which justifies the strictness in the inequality (\ref{2}). $\Box$ 

\begin{remark}\label{Gram}
In different terms, Lemma~\ref{co-ordered} shows that when the families of (pairwise distinct) vectors $\{u_i, 1\le i\le s\}$ and $\{v_i, 1\le i\le s\}$ are co-ordered, $s\times s$ Gram matrix $(\langle u_i, v_j\rangle)_{1\le` i,j\le s}$ is tropically non-singular \cite{MS} with the strict minimum of the sum of diagonal entries among all $s!$ matchings. 
\end{remark}
 


\begin{proposition}\label{dimension}
When $V\subset \RR^n$ is a rational polyhedron, it holds $TH_V(k)=t(k)$.
\end{proposition}  

{\bf Proof}. In the direction that $TH_V(k)\le t(k)$, the statement is clear. \vspace{1mm}

Let us prove the inverse direction. Choose $m$ coordinates among $x_1,\dots,x_n$ such that the projection of $L$ onto the subspace spanned by the chosen coordinates, has the (full) dimension $m$. W.l.o.g. one can suppose that those are just the first coordinates $x_1,\dots,x_m$. Then one can consider  $x_1,\dots,x_m$ as the coordinates on $V$ and pick a tropical monomial in each equivalence class (with respect to $L$). Thus, one can represent any min-plus polynomial from $T_n(k)$ restricted on $V$ as follows:
\begin{equation}\label{1}
\min_{1\le i\le t(k)} \{\sum_{1\le j\le m} b_{i,j}x_j +b_{i,0}\}    
\end{equation}
where $b_{i,0}\in \RR\cup \{\infty\}$ and $b_{i,j}\in \QQ \cup \{\infty\}$ (we draw attention that $b_{i,j}$ can be rationals rather than integers). Observe that $b_{i,j}$ for $1\le j\le m$ are fixed, while $b_{i,0}, 1\le i\le t(k)$ vary.

Our goal is for each $1\le i\le t(k)$ to find a point $v_i=(v_{i,1},\dots,v_{i,m})\in V$ such that the minimum in (\ref{1}) at $v_i$ is (uniquely) attained for $i$, i.e. at the linear function $\sum_{1\le j\le m} b_{i,j}x_j +b_{i,0}$ for suitable $b_{i,0} \in \RR$ (see Definitions~\ref{dependent}, ~\ref{tropical_hilbert}). 

One can pick points $v_i\in V, 1\le i\le t(k)$ such that the families $\{b_i:=(b_{i,1},\dots,b_{i,m})\ :\ 1\le i\le t(k)\}$ (see (\ref{1})) and $\{v_i, 1\le i\le t(k)\}$ are co-ordered, taking into account that $\dim V=m$.

We apply Lemma~\ref{diagonal} to the matrix $(w_{i,j}:=\langle b_i,v_j\rangle)_{1\le i,j\le t(k)}$ (see Lemma~\ref{co-ordered}), and in (\ref{1}) put $b_{i,0}:=w_i,\ 1\le i\le t(k)$. Then in the resulting min-plus polynomial defined by (\ref{1}), the minimum in (\ref{1}) for the point $v_i, 1\le i\le t(k)$ is attained at the linear function $\sum_{1\le j\le m}b_{i,j}x_j+b_{i,0}$. This completes the proof of Proposition~\ref{dimension}. $\Box$

\begin{remark}\label{semigroup}
For a polyhedron $V\subset \RR^n$ since the equivalence classes (with respect to $L$) of $\ZZ_{\ge 0}^n$ form a (finitely generated) semigroup, \cite{Khovanskii} implies that $t(k):=t_L(k)$ coincides with a polynomial in $k$ for sufficiently big $k$. The degree of this polynomial equals $m=\dim V$.
\end{remark}

Now we explicitly compute $t_V(k)$ in case of an interval $V\subset \RR^2$ of a line in the plane. 

\begin{proposition}\label{line}
For an interval $V$ of a line $\{pX_1+qX_2=r\} \subset \RR^2,\ r\in \RR$, where either integers $p,q$ are relativly prime or $pq=0, p+q=1$,  it holds $t_V(k)=(|p|+|q|)(k-1)+1$.


\end{proposition}

{\bf Proof}. Equivalence classes of $M_k$ correspond bijectively to lines of the form $\{qX_1-pX_2=c,\ c\in \ZZ\}$ which have non-empty intersection with $M_k$. If $p,q\ge 0$ then $-p(k-1)\le c\le q(k-1)$. Otherwise, if $q>0, p\le 0$ then $0\le c\le (q-p)(k-1)$. $\Box$



I

\begin{remark}
1. If $V$ is a line then $V$ is a parallel shift of the tropical variety $trop(X_1^pX_2^q+1)$ when $p,q\ge 0$
  or $trop(X_1^{-p}+X_2^q)$ when $p\le 0, q\ge 0$ (cf. Proposition~\ref{line}) . 

2. When $p,q$ are rationally independent, we have $t_V(k)=k(k+1)/2$ since each equivalence class consists of a single point. 
\end{remark}

\begin{remark}\label{Hilbert}
If a min-plus prevariety $V=\cup_l V_l$ where $V_l$ is a polyhedron for each $l$, then for every $l_0$ it holds
$$TH_{V_{l_0}}\le TH_V\le \sum_l TH_{V_l}.$$
Thus, due to Remark~\ref{semigroup} $TH_V$ grows between two polynomials of the degree $\dim V$.
\end{remark}




\begin{proposition}\label{zero}
Let $V\subset \RR^n$ consist of $s$ points. Then $TH_V(k)=s$ for  $k\ge ns$.
\end{proposition}

{\bf Proof}. The proof is similar to the proof of Proposition~\ref{dimension}.

Let $V=\{v_1,\dots,v_s\}$. One can pick a family of vectors $b_1,\dots,b_s\in \ZZ_{\ge 0}^n$ such that this family is co-ordered with $V$ (moreover, one can pick $b_1,\dots,b_s$ with coordinates less than $s$). According to Lemmas~\ref{co-ordered}, ~\ref{diagonal} there exist $b_{i,0}\in \RR, 1\le i\le s$ such that the min-plus polynomial (\ref{1}) (replacing in it $m$ by $n$ and $t(k)$ by $s$, respectively) attains the strict minimum for a point $v_i, 1\le i\le s$ at the tropical monomial $\langle b_i, x\rangle +b_{i,0}$. 
Hence $TH_V(k)\ge s$ (for $k\ge ns$).

On the other hand, $TH_V(k)\le s$ due to Remark~\ref{Hilbert}. $\Box$

\section{Existence of the tropical degree of a one-dimensional tropical prevariety}

For a tropical prevariety $V\subset \RR^n$ with $\dim(V)=1$ we define its {\it tropical degree}  as the limit
$$degT (V):=\lim_{k\to \infty} \frac{TH_V(k)}{k}.$$
\noindent The following theorem justifies the existence of this limit.

\begin{theorem}\label{degree}
For a tropical prevariety $V\subset \RR^n$ with $\dim V=1$ its tropical degree is defined.     
\end{theorem}

{\bf Proof}. For the sake of simplifying notations assume w.l.o.g. that $n=2$.

We need the following extension of the subadditivity lemma \cite{S}.

\begin{lemma}\label{subadditivity}
Let $\{0\le a_i\in \RR, i\ge 1\}$ be a non-decreasing sequence fulfilling a condition
$$a_{kr}\ge r(a_k-c), k,r\ge 1$$
\noindent for a constant $c\in \RR$. Then there exists a limit
$$\lim_{q\to \infty} \frac{a_q}{q}.$$
\end{lemma}

{\bf Proof of the lemma}. Denote 
$$M_0:= \limsup_{q\to \infty}  \frac{a_q}{q}.$$
When $M_0=0$ the lemma is obvious. Let $M_0>0$. Take an arbitrary $M\in \RR,\ 0<M\le M_0$. Fix an arbitrary $0<e<M/3$ for the time being. Take sufficiently big $k$ such that $a_k/k>M-e,\, c/k<e$. One can represent any $q=kp+q_0,\ 0\le q_0<k$. It holds
$$\frac{a_q}{q}\ge \frac{pa_k-cp}{kp+q_0}>\frac{pa_k-cp}{k(p+1)}.$$
\noindent Taking any $p>M/e$, we obtain 
$$\frac{pa_k}{k(p+1)}>\frac{p(M-e)}{p+1}>M-2e,\ \frac{cp}{k(p+1)}<e,$$
\noindent we conclude that $\frac{a_q}{q}>M-3e$. $\Box$ \vspace{2mm}

We continue the proof of the theorem.

Fix an integer $k\ge 1$ for the time being and assume that the tropical monomials 
\begin{eqnarray}\label{90}
\{ix+jy\, |\, 0\le i+j\le k,\, (i,j)\in S\}    
\end{eqnarray}
are independent on $V$ for a suitable set $S\subset \ZZ^2$ such that $|S|=TH_V(k)$. Let
$$f:=\min_{(i,j)\in S} \{ix+jy+a_{ij}\},\ a_{ij}\in \RR$$
be a certificate of independence of (\ref{90}) satisfying Definition~\ref{dependent}.

Slightly perturbing the coefficients $\{a_{ij}\}$ one can suppose w.l.o.g. that at no point $v\in V$ the values of any three linear functions among $\{ix+jy+a_{ij}\, |\, (i,j)\in S\}$ coincide.

Denote by $P\subset \RR^3$ the convex hull of the rays $\{(i,j,a)\, |\, a\ge a_{ij}\}$ for all $(i,j)\in S$. Then $(i,j,a_{ij})$ is a vertex of $P$ for each $(i,j)\in S$ due to the independence of (\ref{90}).

It holds $V=\cup_{1\le l\le c} V_l$ where every $V_l, 1\le l\le c$ is an interval (whatever open or closed, finite or infinite). For every $V_l$ the graph $F_l\subset V_l \times \RR$ of the restriction $f|_{V_l}$ is a convex polygon. We say that two points $(i_1,j_1), (i_2,j_2)\in S$ are {\it adjacent} on $V_l$ if the linear functions
\begin{eqnarray}\label{91}
(i_1x+j_1y+a_{i_1,j_1})|_{V_l},\, (i_2x+j_2y+a_{i_2,j_2})|_{V_l}
\end{eqnarray}
correspond to adjacent edges $E_1, E_2$ of the polygon $F_l$.

\begin{lemma}\label{edge}
If points $(i_1,j_1), (i_2,j_2)\in S$ are adjacent on some $V_l$ then $(i_1,j_1,a_{i_1,j_1}), (i_2,j_2,a_{i_2,j_2})$ are the endpoints of an edge of the polyhedron $P$.      
\end{lemma}

{\bf Proof of the lemma}. Suppose the contrary. Denote by $P'\subset P$ the convex hull of the rays $\{(i,j,a)\, |\, a\ge a_{ij}\}$ for all $(i_1,j_1), (i_2,j_2)\neq (i,j)\in S$. Then the interval with the endpoints $(i_1,j_1,a_{i_1,j_1}), (i_2,j_2,a_{i_2,j_2})$ has a common point 
$$w:=\sum_{(i_1,j_1), (i_2,j_2)\neq (i,j)\in S} \alpha_{ij} (i,j,a_{ij}+b_{ij})\in P'$$
\noindent with $P'$, where $0\le \alpha_{ij},\ \sum_{(i,j)} \alpha_{ij}=1,\, b_{ij}\ge 0$.

Denote by $(v,b)\in F_l,\, v\in V_l$ the (unique) common vertex of $F_l$ of adjacent edges $E_1, E_2$ corresponding to
(\ref{91}), in other words
$$(i_1x+j_1y+a_{i_1,j_1})(v)=(i_2x+j_2y+a_{i_2,j_2})(v)=b.$$
\noindent It holds $(ix+jy+a_{ij})(v)>b$ for each $(i_1,j_1), (i_2,j_2)\neq (i,j)\in S$, which contradicts to that $w$ belongs to the interval with the endpoints $(i_1,j_1,a_{i_1,j_1}), (i_2,j_2,a_{i_2,j_2})$.
The achieved contradiction proves the lemma. $\Box$ \vspace{2mm}

We come back to the proof of the theorem. Introduce a graph $G$ whose vertices are the vertices of the polyhedron $P$. Two vertices 
$(i_1,j_1,a_{i_1,j_1}), (i_2,j_2,a_{i_2,j_2})$
are connected by an edge in $G$ if $(i_1,j_1),\, (i_2,j_2)$ are adjacent on $V_l$ for some $1\le l\le c$. 
Lemma~\ref{edge} implies that $G$ is a subgraph of the graph of edges of the polyhedron $P$.


\begin{lemma}\label{connected}
The graph $G$ has at most $c$ connected components.    
\end{lemma}

{\bf Proof of the lemma}. Take any $c+1$ points $(i_0,j_0),\dots, (i_c,j_c)\in S$. For each $0\le p\le c$ there exists (perhaps, not unique) an edge of $F_l$ for some $1\le l\le c$ given by the linear function
$$(i_px+j_py+a_{i_p,j_p})|_{V_l}.$$
\noindent There exist $(i_p,j_p),\, (i_q,j_q),\, p\neq q$ whose corresponding edges belong to the same $F_l$ for 
appropriate $1\le l\le c$.
Therefore, the vertices $(i_p,j_p,a_{i_p,j_p}),\, (i_q,j_q,a_{i_q,j_q})$ lie in the same connected component of $G$. $\Box$ \vspace{2mm}

We return to the proof of the theorem. Fix an integer $r\ge 1$ for the time being. Our goal is to produce a set 
$$S^{(r)}\subset \{(i,j)\, |\, 0\le i+j\le kr\}\subset \ZZ^2,\ S^{(r)}\supset rS$$
\noindent and a tropical polynomial
$$f^{(r)}:=\min_{(i,j)\in S^{(r)}} \{ix+jy+b_{ij}\}$$
\noindent being a certificate of independence of   the tropical monomials $\{ix+jy\, |\, (i,j)\in S^{(r)}\}$ on $V$ (see Definition~\ref{dependent}).

Consider any pair of vertices $(i_1,j_1,a_{i_1,j_1}), (i_2,j_2,a_{i_2,j_2})$ of $P$ such that the linear functions
$$i_1x+j_1y+a_{i_1,j_1}|_{V_l},\ i_2x+j_2y+a_{i_2,j_2}|_{V_l}$$
\noindent correspond to adjacent edges $E_1, E_2$ of the polygon $F_l$ for some $l$. Denote by $(v,b)\in F_l$ the (unique) vertex being common for the edges $E_1, E_2$. 
It means that there is an edge in the graph $G$ connecting vertices $(i_1,j_1,a_{i_1,j_1}), (i_2,j_2,a_{i_2,j_2})$. In particular, it holds
$$(i_1x+j_1y+a_{i_1,j_1})(v)=(i_2x+j_2y+a_{i_2,j_2})(v)=b.$$

The set $S^{(r)}$ consists of the points
$$\{(r-p)(i_1,j_1)+p(i_2,j_2)\ |\ 0\le p\le r\}$$
for every edge $(i_1,j_1,a_{i_1,j_1}), (i_2,j_2,a_{i_2,j_2})$ of the graph $G$.

Now we define the coefficients $b_{ij}$ of $f^{(r)}$. For each produced point
$$(i,j):=(r-p)(i_1,j_1)+p(i_2,j_2),\ 0\le p\le r$$
\noindent we put
\begin{eqnarray}\label{92}
b_{ij}:=(r-p)a_{i_1,j_1}+pa_{i_2,j_2}+\varepsilon p(r-p)
\end{eqnarray}
\noindent for an arbitrary sufficiently small $\varepsilon >0$. Below we discuss, how small should be $\varepsilon$. In fact, in place of $p(r-p)$ one can take any strictly convex function which vanishes at $p=0, r$. 

Consider any point $(v_1,b_1)$ being a vertex of the polygon $F_m$ for some $m$ such that either $(i_1x+j_1y+a_{i_1,j_1})(v_1)>b_1$ or $(i_2x+j_2y+a_{i_2,j_2})(v_1)>b_1$. It means that either $(i_1x+j_1y+a_{i_1,j_1})|_{V_m}$ or $(i_2x+j_2y+a_{i_2,j_2})|_{V_m}$ does not correspond to an edge of the polygon $F_m$ incident to the vertex $(v_1,b_1)$. Then (\ref{92}) implies that 
$$(ix+jy+b_{ij})(v_1)>rb_1, 1\le p<r$$
\noindent for sufficiently small $\varepsilon$.

It holds $f^{(r)}\le rf$, and $f^{(r)}|_{V}$ coincides with $rf|_V$ out of certain small neighbourhoods of the points $v_1$ such that $(v_1,b_1)$ is a vertex of $F_m$ for some $m$.

Moreover, for sufficiently small $\varepsilon$ in a suitable neighbourhood of $v$ in $V_l$ the function $f^{(r)}$ coincides with (see (\ref{92}))
$f_{(i_1,j_1), (i_2,j_2)}:=$
$$\min_{0\le p\le r} \{((r-p)i_1+pi_2)x+((r-p)j_1+pj_2)y+(r-p)a_{i_1,j_1}+pa_{i_2,j_2}+\varepsilon p(r-p)\}.$$
\noindent The graph $F_{(i_1,j_1), (i_2,j_2)}\subset V_l\times \RR$ of the concave function $f_{(i_1,j_1), (i_2,j_2)}|_{V_l}$ is a convex polygon 

(i) with $r+1$ edges and

(ii) $r$ vertices in a neighbourhood of the point $(v, 
rb)$.

The proof of (i) relies on the following obvious lemma.

\begin{lemma}
The graph of a concave function on $\RR$ given by
$\min_{0\le l\le r} \{a_l+lx\}$, where the sequence $\{a_l\}$ is strictly convex, i.e. $2a_l<a_{l-1}+a_{l+1}, 1\le l<r$, is a polygon with $r$ vertices.
\end{lemma}

To verify (ii) denote $V_l:=\{v+tw\}$ for some vector $w\in \ZZ^2$, where $t\in \RR$ varies in a certain interval (which contains $0$). It holds
$$\langle (i_1,j_1), v\rangle +a_{i_1,j_1}=\langle (i_2,j_2), v\rangle +a_{i_2,j_2}=b,\ \langle (i_1-i_2,j_1-j_2), w\rangle \neq 0.$$
\noindent For $0\le p<r$ the $p$-th vertex $v+t_pw$ of the polygon $F_{(i_1,j_1), (i_2,j_2)}$ satisfies the following equation:
$$\langle ((r-p)i_1+pi_2, (r-p)j_1+pj_2),\ (v+t_pw)\rangle + (r-p)a_{i_1,j_1}+pa_{i_2,j_2}+\varepsilon p(r-p) =$$
$$\langle ((r-p-1)i_1+(p+1)i_2, (r-p-1)j_1+(p+1)j_2),\ (v+t_pw)\rangle$$
$$+ (r-p-1)a_{i_1,j_1}+(p+1)a_{i_2,j_2}+\varepsilon (p+1)(r-p-1).$$
\noindent Hence 
$$t_p=\frac{\varepsilon (r-2p-1)}{\langle (i_1-i_2, j_1-j_2), w\rangle}.$$
\noindent Therefore, for sufficiently small $\varepsilon$ the point $v+t_pw\in V_l$ is located in a small neighbourhood of $v$.

Thus, the tropical polynomial $f^{(r)}$ is a certificate of independence of tropical monomials $\{ix+jy\ |\ (i,j)\in S^{(r)}\}$ on $V$ (cf. Definition~\ref{dependent}).

Due to Lemma~\ref{connected} it holds
$$TH_V(kr)\ge |S^{(r)}|\ge (|S|-c)r=(TH_V(k)-c)r.$$
\noindent Lemma~\ref{subadditivity} completes the proof of the Theorem~\ref{degree}. $\Box$

\begin{remark}
If to stick with the conjecture that $TH_V(k)$ is (eventually) a polynomial then  we obtain that $degT (V)$ is the leading coefficient of the linear (eventual) polynomial $TH_V(k)$.
\end{remark}

\section{Tropical degree of a star}


In this section for the sake of simplyfing calculations we consider a slight modification of the tropical Hilbert function. Namely, we define $TH_V^{(\Box)}(k)$ to be the maximal number of independent on $V$ tropical monomials among
$$\{i_1x_1+\cdots +i_nx_n\ |\ 0\le i_j\le k, 1\le j\le n\}.$$
\noindent For $\dim(V)=1$ the proof of Theorem~\ref{degree} goes literally through and justifies the existence of the limit
$$degT^{(\Box)}(V):= \lim_{k\to \infty} \frac{TH_V^{(\Box)}(k)}{k}.$$

\begin{definition}\label{star}
Let $V\subset \RR^2$ be a {\it star}, i.e. the union of rays $V_l, 1\le l\le M$ with a common apex.
Let $(p_l, q_l)$ be the spanning vector of $V_l$. We assume that for each $1\le l\le M$ it holds that either $p_l, q_l$  are relatively prime integers or $p_lq_l=0, |p_l|+|q_l|=1$.
    \end{definition}

\begin{theorem}\label{leading}
For a star $V$ its leading coefficient does exist and equals
\begin{eqnarray}\label{50}
degT^{(\Box)} (V)=\max\{\sum_{p_l>0} p_l,\ -\sum_{p_l<0} p_l\} + \max\{\sum_{q_l>0} q_l,\ -\sum_{q_l<0} q_l\} =:B.    \end{eqnarray}
\end{theorem}

{\bf Proof}. Nevertheless, one can assume that the common apex of the rays of $V$ is the origin $(0,0)$. 

{\bf Upper bound}.

First we establish an upper bound on $TH_V^{(\Box)}(k)$. Consider a tropical polynomial 
$$f:= \min_{0\le i,j\le k} \{ix+jy+a_{ij}\}, a_{ij}\in \RR\cup \{0\}.$$
\noindent Assume that $a_{i_0, j_0}=\min_{0\le i,j\le k} \{a_{ij}\}$. Observe that each restriction $f|_{V_l}, 1\le l\le M$ is a concave piece-wise linear function.

The number of slopes of the restriction $f|_{V_l}$ does not exceed

$\bullet$ $(i_0-k)p_l+(j_0-k)q_l$ if $p_l, q_l\le 0$ (cf. Proposition~\ref{line}) since the slope of the restriction $f|_{V_l}$ at $0$ is greater or equal to $-i_0p_l-j_0q_l$;

$\bullet$ $i_0p_l+j_0q_l$ if $p_l, q_l\ge 0$;

$\bullet$ $i_0p_l+(j_0-k)q_l$ if $p_l\ge 0, q_l\le 0$;

$\bullet$ $(i_0-k)p_l +j_0q_l$ if $p_l\le 0, q_l\ge 0$.

Thus, we get inequalities 
\begin{eqnarray}\label{52}
 TH_V^{(\Box)}(k)-1\le \max_{0\le i,j\le k} \{\sum_{p_l> 0} ip_l + \sum_{p_l< 0} (i-k)p_l +\sum_{q_l> 0} jq_l + \sum_{q_l< 0} (j-k)q_l\}=kB.
\end{eqnarray}

{\bf Lower bound}.

W.l.o.g. assume that 
$$-\sum_{p_l <0} p_l\ge \sum_{p_l>0} p_l,\ \ -\sum_{q_l <0} q_l\ge \sum_{q_l>0} q_l.$$
\noindent Note that the maximum of the right-hand side of (\ref{52}) being a linear function in $i,j$, is attained at one of 4 points $(i,j)=(0,0), (0,k), (k,0), (k,k)$. The latter assumption means that the maximum is attained at the point $(i,j)=(0,0)$. 
Other 3 cases one can consider in a similar way.

Suppose that among vectors $(p_l, q_l), 1\le l\le M$ there are $m\le M$ vectors such that either $p_l<0$ or $q_l<0$. Denote those vectors by $(u_1, v_1),\dots, (u_m,v_m)$. Thus, $B=-\sum_{u_l<0} u_l -\sum_{v_l<0} v_l$ (see (\ref{50})).

For integers $c_1,\dots, c_m$ consider a convex polygon
$$P:=P_k(c_1,\dots, c_m):=\{u_lx+v_ly\ge c_l\ |\ 1\le l\le m\} \cap Q_k \subset \RR^2,$$
\noindent where $Q_k:=\{0\le x,y\le k\}$ denotes the square. The edges of $P$ are of 2 types: either located on the edges of $Q_k$ (we call them static) or located on the lines of the form $\{u_lx+v_ly=c_l\}$ (we call them dynamic). Points of an edge being not vertices, we call inner points of this edge.

\begin{lemma}\label{edges}
i) There exists an integer $C:=C(V)$ such that under an assumption that no dynamic edge of the polygon $P$ contains an inner integer point, the sum of lengths of dynamic edges is less than $C$.

ii) Moreover, under the same assumption when $k>C$ one of 2 cases can emerge. In the first case the dynamic edges form a single interval $e_1 \dots e_s$ (on the boundary of $P$ in the clock-wise order), where the left vertex of $e_1$ lies on the static edge $\{x=0, 0\le y\le k\}$, while the right vertex of $e_s$ lies on the static edge $\{0\le x\le k, y=0\}$. We say in this case that the interval $e_1 \dots e_s$ connects the static edges $\{x=0, 0\le y\le k\}$ and $\{0\le x\le k, y=0\}$.

In the second case the dynamic edges form at most of 3 intervals such that

$\bullet$ the first interval connects the static edges $\{x=0, 0\le y\le k\}$ and $\{0\le x\le k, y=k\}$, and each edge of the interval lies on a line $\{u_lx+v_ly=c_l\}$ with $u_l>0, v_l<0$;

$\bullet$ the second interval connects the static edges $\{0\le x\le k, y=k\}$ and $\{x=k, 0\le y\le k\}$, and each edge of the interval lies on a line $\{u_lx+v_ly=c_l\}$ with $u_l<0, v_l<0$;

$\bullet$ the third interval connects the static edges $\{x=k, 0\le y\le k\}$ and $\{0\le x\le k, y=0\}$, and each edge of the interval lies on a line $\{u_lx+v_ly=c_l\}$ with $u_l<0, v_l>0$.
\end{lemma}

{\bf Proof}. The assertion i) of the lemma follows from the observation that any interval with the length greater than $\sqrt{u_l^2+v_l^2}$ of a line $\{u_lx+v_ly=c\}$ for an integer $c$, contains an integer point. Thus, one can take $C(V):=\lceil \sum_{1\le l\le m} \sqrt{u_l^2+v_l^2}\rceil$.

The assertion ii) one verify straight-forwardly. $\Box$ \vspace{2mm}

Now we describe a recursive construction of a set of points $W\subset Q_k \cap \ZZ^2$ and of integers $c_1,\dots,c_m$. As a base of recursion we take $c_l:= C\min\{u_l, v_l\}, 1\le l\le m$ and $W=\emptyset$.

For the recursive step assume that the current integers $c_1,\dots,c_m$ and a set $W$ are already constructed. 
The set $W$ is contained strictly inside the polygon $P:=P_k(c_1,\dots,c_m)$.
Choose a dynamic edge of $P:=P_k(c_1,\dots,c_m)$ lying on a line $\{u_lx+v_ly=c_l\}$ such that this edge contains an integer inner point $w$, provided that such an edge does exist (if there is no such an edge, the recursive construction halts). Then for the next recursive step we put $(c_1,\dots,c_m):= (c_1,\dots, c_{l-1}, c_l-1, c_{l+1},\dots, c_m)$ and $W:=W\sqcup \{w\}$. We say that the line $\{u_lx+v_ly=c_l\}$ is {\it active} at this recursive step. Obviously $P_k(c_1,\dots,c_m) \subset P_k(c_1,\dots, c_{l-1}, c_l-1, c_{l+1},\dots, c_m)$.

Now assume that a recursive step can't be carried out, so the recursion halts. Then  the first case in Lemma~\ref{edges}~ii) can't take place because of the choice of $c_1,\dots,c_m$ at the base of recursion. According to the second case in Lemma~\ref{edges}~ii) one can assume that the current polygon $P$ contains an interval of dynamic edges $e_1 \dots e_s$ (perhaps, empty) such that the left vertex of $e_1$ lies on the static edge $\{x=0, 0\le y\le k\}$, while the right vertex of $e_s$ lies on the static edge $\{0\le x\le k, y=k\}$. Thus, for every $1\le l\le m$ such that $u_l\ge 0, v_l<0$ the line $\{u_lx+v_ly=c_l\}$ intersects the static edge $\{0\le x\le k, y=k\}$ at a point $(x_0, k)$ for which it holds $x_0<C$. Hence $c_l<u_lC+v_lk$. 

The other 2 possible intervals one can study in a similar manner. 

Thus, we obtain the following properties of the described recursive construction.

\begin{corollary}\label{set}
For every $1\le l\le m$ the line $\{u_lx+v_ly=c\}$ is active at some step of the described recursive construction for:

$\bullet$  $v_lC\ge c\ge u_lC+v_lk$ when $u_l\ge 0, v_l<0$;

$\bullet$ $(u_l+v_l)C\ge c\ge (u_l+v_l)(k-C)$ when $u_l, v_l<0$;

$\bullet$ $u_lC\ge c\ge u_lk+v_lC$ when $u_l<0, v_l\ge 0$.

Thus, it holds
$$|W|\ge -k(\sum_{u_l<0} u_l+ \sum_{v_l<0} v_l)+C(\sum_{u_l\ge 0, v_l<0} (v_l-u_l)+ 2\sum_{u_l,v_l<0} (u_l+v_l) +\sum_{u_l<0, v_l\ge 0} (u_l-v_l)).$$
Denoting the right-hand side of the latter inequality by $kB-D$ (see (\ref{50})), where $D:=D(V)$, we conclude that $|W|\ge kB-D$.
\end{corollary}

For every point $w=(w^{(1)}, w^{(2)})\in W$ consider a tropical monomial $w^{(1)}x+w^{(2)}y$.

\begin{lemma}\label{monomials}
The tropical monomials $\{w^{(1)}x+w^{(2)}y\}, w\in W$ are independent on the star $V$.
\end{lemma}

{\bf Proof}. Denote by $w_s\in W, 1\le s\le |W|$ the point which is added to $W$ at the $s$-th step of the described above recursive construction. Consider a tropical polynomial
$$g:=\min_{1\le s\le |W|} \{w^{(1)}_s x+ w^{(2)}_s y + 2^s\}.$$

Fix $1\le l\le m$ for the time being. Denote by $V_l$ the ray spanned by the vector $(u_l, v_l)$,  by $z$ the coordinate of the ray $V_l$ and by $1\le s_c\le |W|$ the number of the recursive step at which the line $\{u_lx+v_ly=c\}$ is active (provided that $s_c$ is defined). Then the restriction
\begin{eqnarray}\label{51}
g|_{V_l}=\min_{0\ge \gamma > c_1\ge c\ge c_2} \{cz+2^{s_c}, \gamma z +\alpha_{\gamma}\}
\end{eqnarray}
\noindent (cf. Proposition~\ref{line}) for integers $\gamma$ and some $0<\alpha_{\gamma} \in \RR\cup\{\infty\}$. The bounds $c_1,c_2$ are taken from Corollary~\ref{set}, i.e.

$\bullet$ $c_1=v_lC, c_2=u_lC+v_lk$ in case $u_l\ge 0, v_l<0$;

$\bullet$ $c_1=(u_l+v_l)C, c_2=(u_l+v_l)(k-C)$ in case $u_l, v_l <0$;

$\bullet$ $c_1=u_lC, c_2=u_lk+v_lC$ in case $u_l<0, v_l\ge 0$.

\begin{lemma}\label{concave}
Let $h=\min_{0\le i\le n} \{iz+a_i\}$ be a tropical univariate polynomial on the non-negative ray $\RR_{\ge 0}$. Assume that $0<a_j\in \RR\cup \{\infty\}, t+1\le j\le n$ and $0<a_i\in \RR, 0\le i\le t, a_{i-1}\ge 2a_i, 1\le i\le t$. Then the slopes $iz+a_i, 0\le i\le t$ occur in $h$. 

Moreover, for a suitable point $0\le z^{(0)}< a_{t-1}-a_t$ it holds 
$$h|_{\RR_{\ge z^{(0)}}} =\min_{0\le i\le t} \{iz+a_i\}|_{\RR_{\ge z^{(0)}}},$$
and all the slopes $iz+a_i, 0\le i\le t$ occur in the piece-wise linear function $h|_{\RR_{\ge z^{(0)}}}$.
\end{lemma}

{\bf Proof}. For each $1\le i\le t$ denote $z_i:=a_{i-1}-a_i>0$. Obviously, $iz_i+a_i=(i-1)z_i+a_{i-1}$. We claim that $jz_i+a_j>iz_i+a_i$ for all $0\le j\le n, j\neq i, i-1$ (i.e. $h(z_i)=iz_i+a_i$). If $j\ge i+1$ then $(j-i)z_i\ge a_{i-1}-a_i>a_i-a_j$. Otherwise, if $j\le i-2$ then 
$$a_j\ge a_{i-1}2^{i-1-j}\ge a_{i-1}(i-j)>(a_{i-1}-a_i)(i-j)+a_i.$$

For $j\ge t+1$ denote by $z_j^{(0)}$ the point such that $jz_j^{(0)}+a_j=tz_j^{(0)}+a_t$. Then it holds $z_j^{(0)}\le a_t-a_j<a_{t-1}-a_t$. One can put $z^{(0)}:=\max_{t+1\le j\le n} \{z_j^{(0)}, 0\}$. $\Box$ \vspace{2mm}

Applying Lemma~\ref{concave} to the tropical polynomial $g|_{V_l}$ (see (\ref{51})) we conclude that the piece-wise linear function $g|_{V_l}$ has at least $c_1-c_2+1$ slopes which correspond to the lines $\{u_lx+v_ly=c\}, c_1\ge c\ge c_2$, which in their turn correspond to the tropical monomials $w_{s_c}^{(1)}x+w_{s_c}^{(2)}y$. This completes the proof of Lemma~\ref{monomials}. $\Box$ \vspace{2mm}

Therefore Corollary~\ref{set} implies the inequality $TH_V^{(\Box)}(k)\ge kB-D$, hence $degT^{(\Box)} (V)= B$ taking into account (\ref{52}), which completes the proof of   Theorem~\ref{leading}. $\Box$


\begin{remark}
One can extend Theorem~\ref{star} to a star $V\subset \RR^n$ being a union of rays spanned by vectors
$$\{(p_{j1},\dots,p_{jn})\in \ZZ^n\ |\ 1\le j\le m\}$$
\noindent where $p_{j1},\dots,p_{jn}$ have no nontrivial common divisor, $1\le j\le m$. One can show that
$$degT^{(\Box)} (V)= \sum_{1\le i\le n} \max\{\sum_{j\ :\ p_{ji}>0} p_{ji},\ -\sum_{j\ :\ p_{ji}<0} p_{ji}\}.$$
\end{remark}

\begin{remark}
Theorem~\ref{leading} differs from the picture in the classical algebraic geometry in which the degree of the union of two varieties of  equal dimensions in a generic position is the sum of their degrees.    
\end{remark} \vspace{2mm}

{\bf Acknowledgements}. The author is grateful to C.~G.~L\'opez for interesting discussions.

\end{document}